\documentclass[aos,MSNbibl,seceqn,dvips]{arximspdf}

\doi{10.1214/12-AOS1081} 
\volume{41}
\issue{2}
\pubyear{2013}
\firstpage{464}
\lastpage{483}

\makeatletter
\newcommand{\rrvert}{\vert}
\newcommand{\llvert}{\vert}
\def\cal{\mathcal}
\newcommand{\R}{\mathbb{R}}
\newcommand{\E}{\mathbb{E}}
\renewcommand{\P}{\mathbb{P}}
\newcommand{\N}{\mathbb{N}}
\newcommand{\Var}{\operatorname{Var}}
\newcommand{\Cov}{\operatorname{Cov}}
\newcommand{\trace}{\operatorname{trace}}

\newtheorem{theorem}{Theorem}[section]

\newtheorem{proposition}[theorem]{Proposition}

\newcommand{\eqref}[1]{(\ref{#1})}
\makeatother

\begin{document}
\begin{frontmatter}

\title{On the conditional distributions of~low-dimensional projections
from high-dimensional data}
\runtitle{Conditional distributions of low-dimensional projections}

\begin{aug}
\author[A]{\fnms{Hannes} \snm{Leeb}\corref{}\ead[label=e1]{hannes.leeb@univie.ac.at}}
\runauthor{H. Leeb}
\affiliation{University of Vienna}
\address[A]{Department of Statistics\\
University of Vienna\\
Universit\"{a}tsstr. 3/5\\
Vienna, 1010\\
Austria\\
\printead{e1}} 
\end{aug}

\received{\smonth{10} \syear{2012}}

%
\begin{abstract}
We study the conditional distribution of low-dimensional projections from
high-dimensional data, where the conditioning is on other low-dimensional
projections. To fix ideas, consider a random
$d$-vector $Z$ that has a Lebesgue density and that is standardized
so that $\E Z =0$ and $\E Z Z' = I_d$.
Moreover, consider two projections defined by unit-vectors $\alpha$ and
$\beta$, namely a response $y = \alpha' Z$
and an explanatory variable $x = \beta' Z$.
It has long been known that
the conditional mean of $y$ given $x$
is approximately linear in $x$, under some regularity conditions;
cf. Hall and Li [\textit{Ann. Statist.} \textbf{21} (1993) 867--889].
However, a corresponding result for the conditional variance has not been
available so far.
We here show that the conditional variance of $y$ given $x$
is approximately constant in $x$ (again, under some regularity conditions).
These results hold uniformly in $\alpha$ and for most $\beta$'s, provided
only that the dimension of $Z$ is large.
In that sense, we see that most linear submodels
of a high-dimensional overall model are approximately correct.
Our findings provide new insights in a variety of modeling scenarios.
We discuss several examples, including
sliced inverse regression,
sliced average variance estimation,
generalized linear models under potential link violation,
and
sparse linear modeling.
\end{abstract}

%


\end{frontmatter}

\section{Introduction}
\label{sec1}

\subsection{Informal summary}

We analyze a situation where a simple model is used
when the true model is, in fact, much more complex.
This situation is particularly common with many contemporary datasets
where the
number of potentially important covariates
or the number of parameters exceeds the sample size;
examples in, say, genomics or economics abound.
When facing a large number of potentially important covariates
or parameters, and a small sample size,
the search for simple models is typically motivated by either one
of two types of assumptions: Parametric assumptions, which postulate
that the true data-generating process is given by
a simple finite-dimensional model;
and nonparametric assumptions, which postulate that the
true data-generating process can be approximated, with arbitrary accuracy,
by comparatively simple finite-dimensional models.
In either case, the underlying postulates can be difficult to verify
in practice. And this difficulty
is often further compounded by the relatively small sample size.
The results that we obtain here
provide an alternative justification
for the use of simple models. We analyze a scenario
where most simple submodels are approximately correct,
provided only that the overall model is sufficiently complex,
irrespective of whether or not the true data-generating process
is given by, or can be closely approximated by,
a finite-dimensional simple model.
This is the main conceptual contribution of this paper.

On a technical level, we extend and refine a method pioneered by Hall
and Li~\cite{Hal93a}.
In that reference, the authors propose a novel approach for studying
conditional means of linear projections under
weak distributional assumptions;
cf. Theorem~3.2 of~\cite{Hal93a}.
Our technical core contribution is an extension of Hall and Li's approach
to also cover conditional variances and higher conditional moments,
and a more explicit control of error terms that allows us to prove
strong statements like~\eqref{c1000} and~\eqref{c2000}, which follow.
Note that we deal here with
the largest singular value of conditional covariance matrices
of increasing dimension,
which turns out to be considerably more challenging
than handling the norm of the conditional mean
vectors that are treated in~\cite{Hal93a}.

The paper is organized as follows: We continue this section
with a more detailed overview of our findings, and with a
discussion of some interesting consequences.
In Section~\ref{S4}, we present our main result, namely Theorem~\ref{t0},
and give an outline of its proof.
The proof is based on five basic steps that correspond to
five propositions that are also given in Section~\ref{S4}. The (more technical)
proofs of these propositions are relegated to the supplementary material~\cite{Lee12b}.

\subsection{Overview of results}
\label{S2}

Consider a random $d$-vector $Z$ that has a Lebes\-gue density and that is
standardized so that
$\E Z = 0$ and $\E Z Z' = I_d$.
Throughout, we will study projections of $Z$ of the form
$\alpha' Z$ and $\beta' Z$ for unit $d$-vectors $\alpha$ and $\beta$.
The conditional mean of $\alpha'Z$ given $\beta'Z=x$ will be denoted
by $\E[\alpha'Z\|\beta'Z = x]$; other conditional expectations are
defined similarly.\footnote{\label{fmeas}There is a measurable function $g\dvtx \R\to\R$
so that $\E[ \alpha' Z \| \beta' Z] = g(\beta' Z)$ holds,
and we write $\E[ \alpha' Z \| \beta' Z = x]$ for $g(x)$;
the existence of $g$ is guaranteed by, say,
Theorem~4.2.8 of~\cite{Dud02a}.
When we say that $\E[ \alpha' Z \| \beta' Z = x]$ is linear
in $x$ [as in condition~(i), which follows], we mean that $g(x)$
can be chosen to be linear.
Similar considerations apply, mutatis mutandis, to
expressions like $\E[Z\|\beta'Z=x]$,
$\E[ (\alpha' Z)^2 \|\beta'Z = x]$ or
$\E[ Z Z' \|\beta'Z = x]$.}
Our main results are concerned with the conditional mean
and with the conditional variance of $\alpha' Z$ given that
$\beta' Z = x$.
To introduce these,
consider the following two conditions: The vector $\beta$ is such that:

\begin{longlist}[(ii)]
\item[(i)] 
for each $\alpha$,
the conditional mean of $\alpha'Z$ given $\beta' Z = x$
is linear in $x\in\R$;
\item[(ii)] 
for each $\alpha$,
the conditional variance of $\alpha' Z$ given $\beta' Z = x$
is constant in $x\in\R$.
\end{longlist}

Suppose that $\alpha$ is an unknown parameter and that
one can observe $y$ and $x$ given by $y = \alpha'Z$ and $x=\beta'Z$,
respectively.
If $\beta$ is such that both (i) and (ii) hold, then $y$ can be
decomposed into
the sum of a linear function of $x$ and a remainder-term, or error-term,
whose conditional mean given $x$ is zero and whose conditional variance
given $x$ is constant; in other words, the model
\[
y = \gamma x + u
\]
applies,
where $\E[u\|x] = 0$ and $\Var[u\|x]=\sigma^2$, and
where $\gamma\in\R$ and $\sigma^2\geq0$ are unknown parameters
that are given by
$\gamma=\alpha'\beta$ and $\sigma^2 = 1-(\alpha'\beta)^2$, respectively.
(Indeed, (i) implies that $\E[y\|x] = \E[\alpha'Z\|\beta'Z] =
\mu+\gamma(\beta'Z)$ for some real constants $\mu$ and $\gamma$.
It now follows from $\E Z = 0$ that
$\mu= 0$, and $\E Z Z' = I_d$ implies that $\gamma= \alpha'\beta$.
Moreover, (ii) entails that $\Var[\alpha'Z\|\beta'Z] = \sigma^2$
for some constant~$\sigma^2$; hence
$1 = \Var(y) = \Var(\alpha'Z) = \E[ \Var[\alpha'Z \| \beta'Z]] +
\E[ (\E[\alpha'Z\|\beta'Z])^2] = \sigma^2 + (\alpha'\beta)^2$, so
that $\sigma^2$ is given by $\sigma^2 = 1-(\alpha'\beta)^2$.)
These observations continue to hold also if the vectors $\alpha$ and
$\beta$ are not normalized to unit length, mutatis mutandis.

Conditions (i) and (ii) are satisfied for each $\beta$
if $Z$ is normally distributed, that is, $Z\sim N(0,I_d)$. But besides
the Gaussian law, the class of distributions that satisfy~(i) and (ii)
for each $\beta$ appears to be quite small:
Indeed, if $Z$ satisfies (i)
for each $\beta$, then the law of $Z$ is spherically
symmetric~\cite{Eat86a}.
And if, in addition, also (ii) holds for some $\beta$,
then $Z$ is Gaussian~\cite{Bry95a}, Theorem~4.1.4.

Under comparatively mild
conditions on the distribution of $Z$, we here show that
both conditions (i) and (ii) are \emph{approximately} satisfied for
\emph{most} unit-vectors~$\beta$, namely
for a set of unit-vectors $\beta$ whose size, as measured by the uniform
distribution on the unit-sphere in $\R^d$, goes to one as $d\to\infty$.
To state this more formally, we first describe two preliminary results,
namely \eqref{c100} and \eqref{c200}, which follow,
and then extend these to our main results in \eqref{c1000} and \eqref{c2000}
below.

To introduce the two preliminary results mentioned earlier,
we note that (i) and (ii) together
are equivalent to the requirement that
$\E[ Z \| \beta' Z = x]= \beta x$ and
$\E[ Z Z' \|\beta'Z = x] =I_d + (x^2-1) \beta\beta'$ hold for each
$x\in\R$.\vadjust{\goodbreak}
(In other words, the first two moments of the conditional distribution
coincide with what they would be in the Gaussian case.)
From this, it is easy to see that (i) and (ii) are also
equivalent to the requirement that both
%
\begin{eqnarray}
\label{c10} \bigl\Vert \E\bigl[ Z \| \beta' Z = x\bigr] \bigr\Vert^2 -
x^2 &=& 0 \quad\mbox{and}
\\
\label{c20} \bigl\Vert \E\bigl[ Z Z' \|\beta'Z = x\bigr] -
\bigl(I_d + \bigl(x^2-1\bigr) \beta\beta'
\bigr)\bigr \Vert & =& 0
\end{eqnarray}
hold for each $x\in\R$.
Note that we use the notation $\Vert \cdot\Vert $ to denote both
the Euclidean norm of vectors as in \eqref{c10}
and the operator norm of matrices as in \eqref{c20}.
The left-hand side of \eqref{c10} can be written as
\[
\bigl\Vert \E\bigl[ Z \| \beta' Z = x\bigr] \bigr\Vert^2 -
x^2 = \sup_{\alpha} \bigl| \E\bigl[\alpha' Z\|
\beta'Z = x\bigr] - \alpha' \beta x \bigr|^2,
\]
with the supremum taken over all unit-vectors $\alpha\in\R^d$,
as is elementary to verify.\footnote{\label{footA}
This easily follows from the fact that $\E[Z\|\beta'Z]$
(resp., $\beta\beta'Z$)
is the orthogonal projection of~$Z$ into the space of all
measurable (resp., linear) functions of $\beta'Z$ in $L_2(\P)$, and
from the relation between the unconditional and the conditional
variance.}
Hence, the left-hand side of \eqref{c10}
is always nonnegative and
can be interpreted as the worst-case deviation of the regression function
$\E[\alpha' Z\|\beta'Z=x]$ from a linear function
at $x$. The left-hand side of \eqref{c20} can
be interpreted in a similar fashion.
For fixed $x\in\R$, the condition
\eqref{c10} is approximately satisfied for most $\beta$'s if $d$ is large,
under the assumptions of Theorem~3.2 in~\cite{Hal93a} [see also
equation (1.5) in that reference], in the sense that
%
\begin{equation}
\label{c100} \upsilon \bigl\{ \beta\in\R^d\dvtx \bigl\Vert \E\bigl[ Z \|
\beta' Z = x\bigr] \bigr\Vert^2 - x^2 > \varepsilon
\bigr\} \stackrel{d\to\infty} {\longrightarrow} 0
\end{equation}
for each fixed $x\in\R$ and
for each $\varepsilon> 0$, where $\upsilon$ denotes the
uniform distribution on the unit-sphere in $\R^d$.
We here show that condition \eqref{c20} is similarly approximately
satisfied for most $\beta$'s if $d$ is large, in the sense that
%
\begin{equation}
\label{c200} \upsilon \bigl\{ \beta\in\R^d\dvtx\bigl \Vert \E\bigl[ Z
Z' \|\beta'Z = x\bigr] - \bigl(I_d +
\bigl(x^2-1\bigr) \beta\beta'\bigr) \bigr\Vert > \varepsilon
\bigr\} \stackrel{d\to\infty} {\longrightarrow} 0
\end{equation}
for each $x\in\R$ and
for each $\varepsilon> 0$ under the assumptions of Theorem~\ref{t0}
in Section~\ref{S4}.

So far, we have seen for \emph{fixed} $x\in\R$ and for most $\beta$'s
that \eqref{c10} and \eqref{c20} are
approximately satisfied, in the sense that \eqref{c100} and \eqref{c200}
hold under some conditions.
Our main result is that
\eqref{c10} and \eqref{c20} are approximately satisfied
for most $\beta$'s and for most $x$'s:
Under the assumptions of Theorem~\ref{t0},
there are Borel subsets $B_d$ of $\R^d$ satisfying
$\upsilon(B_d) \stackrel{d\to\infty}{\longrightarrow} 1$ so that
%
\begin{eqnarray}
\label{c1000} \sup_{\beta\in B_d} \P \bigl( \bigl\Vert \E\bigl[ Z \| \beta' Z
\bigr] \bigr\Vert^2 - \bigl(\beta'Z\bigr)^2 >
\varepsilon \bigr) &\stackrel{d\to\infty} {\longrightarrow}& 0 \quad \mbox{and}
\\
\label{c2000} \sup_{\beta\in B_d} \P \bigl( \bigl\Vert \E\bigl[ Z Z' \|
\beta'Z\bigr] - \bigl(I_d + \bigl(\bigl(
\beta'Z\bigr)^2-1\bigr) \beta\beta'\bigr) \bigr\Vert > \varepsilon \bigr) &\stackrel{d\to\infty} {\longrightarrow}& 0
\end{eqnarray}
hold for each $\varepsilon> 0$.\vadjust{\goodbreak}

Following a referee's suggestion, we now compare our findings
to the work of Diaconis and Freedman~\cite{Dia84a}, which is an
important precursor to the results of~\cite{Hal93a} and hence, a fortiori,
also to the results in this paper; see also the discussion surrounding
the displays (1.7)--(1.8) in~\cite{Hal93a}.
(Moreover, the recent work of D\"umbgen and Zerial~\cite{Due11a}
should be mentioned here, where several extensions and generalizations of
the results of~\cite{Dia84a} are provided.)
Under the assumptions of Theorem~\ref{t0},
Proposition~5.2 of~\cite{Dia84a} entails, for large $d$, that
the (bivariate) joint distribution of $\alpha'Z$ and $\beta'Z$ is
approximately
normal, with zero means, unit variances, and covariance
$\alpha'\beta$, for most pairs of unit-vectors $\alpha$ and $\beta$
(in the sense of weak convergence in probability with respect to the
product measure $\upsilon\otimes\upsilon$ as $d\to\infty$).
Because the normal distribution has linear conditional means
and constant conditional variances, this suggests, but does not prove,
that
%
\begin{equation}
\label{c150} \P \bigl(\bigl | \E\bigl[ \alpha' Z \| \beta' Z
\bigr] - \alpha'\beta \beta'Z \bigr|^2 >
\varepsilon \bigr) \stackrel{d\to\infty} {\longrightarrow} 0
\end{equation}
for each $\varepsilon>0$ and for most pairs of unit-vectors $\alpha$ and
$\beta$
in $\R^d$ (in the sense of convergence
in probability as a function of $\alpha$ and $\beta$
with respect to $\upsilon\otimes\upsilon$).
If~$\alpha$ is treated as an unknown parameter,
and if the observations $\alpha'Z$ and $\beta'Z$ are treated
as response and explanatory variable, respectively, then \eqref{c150}
entails, for \emph{most $\alpha$'s and $\beta$'s}, that
the response can be approximated by a linear function of
the explanatory variable plus an error term with zero mean conditional
on $\beta'Z$, provided that $d$ is large.
The approximating linear function is $\alpha'\beta(\beta'Z)$.
But for large $d$, we also have $\alpha'\beta\approx0$ for most
$\alpha$'s and $\beta$'s (with respect to $\upsilon\otimes\upsilon$).
The statement in~\eqref{c1000}, on the other hand, is equivalent to
%
\begin{equation}
\label{c151} \P \Bigl( \sup_{\alpha} \bigl| \E\bigl[ \alpha' Z \|
\beta' Z\bigr] - \alpha'\beta \beta'Z
\bigr|^2 > \varepsilon \Bigr) \stackrel{d\to\infty} {\longrightarrow} 0
\end{equation}
for each $\varepsilon>0$ (in probability as a function of $\beta$ with
respect to $\upsilon$); cf. the discussion following \eqref{c20}.
The statement in \eqref{c151} is obviously much
stronger than that in~\eqref{c150}. And it guarantees that
the conditional mean of $\alpha'Z$ is approximately linear in $\beta'Z$,
\emph{for all $\alpha$'s and for most $\beta$'s};
this includes, in particular, the statistically interesting case
where $\alpha$ is parallel, or close to parallel, to $\beta$.
Finally, as already observed
in~\cite{Hal93a},
{``Diaconis and Freedman's result does not provide clues as to whether
} [statements like (1.8)] { might be true or false.}''
Similar observations also apply to conditional variances,
mutatis mutandis.

\subsection{Discussion}

If the left-hand sides of \eqref{c1000} and \eqref{c2000} are both small,
and if $\beta\in B_d$, then the simple linear model, where
the response $\alpha' Z$ is explained by a linear function of
the explanatory variable $\beta' Z$ plus an error that has
zero mean and constant variance given $\beta' Z$, is approximately
correct, irrespective of the unit-vector $\alpha$.
Here, ``approximately correct'' means
that the expressions on the left-hand sides of \eqref{c10} and \eqref{c20}
are at most $\varepsilon$
for a range of values $x$ that contains the explanatory variable $\beta' Z$
with high probability. Under the conditions of Theorem~\ref{t0},
a sufficiently large dimension is enough to guarantee that $B_d$ is large
and that the left-hand sides of \eqref{c1000}
and \eqref{c2000} are small.\footnote{Note that this disentangles the issue of
(approximate) model validity and the issue of model performance:
The model is approximately valid if $\beta\in B_d$,
irrespective of $\alpha$;
the performance of this model, on the other hand,
that is, the performance of $\beta'Z$ as a predictor for $\alpha'Z$,
depends on both $\beta$ and $\alpha$.
Under classical parametric or nonparametric assumptions,
a simple model that is (approximately) correct
typically also performs well.}

The statistical impact of our results is most
pronounced in situations where the sample size is small and the
dimension is large.
Assume that Theorem~\ref{t0} applies, and
consider a collection of $n$ independent copies of the
pair $(\alpha' Z, \beta'Z)$ that we denote by
$(\alpha' Z_i, \beta' Z_i)$, $i=1,\ldots, n$, with $\beta\in B_d$.
If $d$ is large and $n$ is comparatively small, so that
the left-hand sides of both \eqref{c1000} and \eqref{c2000} are
still small even when multiplied by $n$,
then the simple linear model discussed in the preceding paragraph
can also be used to approximately describe
the relation between $\alpha' Z_i$ and $\beta' Z_i$ for each
$i=1,\ldots
, n$,
irrespective of $\alpha$.

We stress that additional data may give reason to dismiss
the simple linear model considered in the preceding paragraphs in favor
of a more complex one,
because the error suffered from using a model that is only approximately
correct will typically become apparent
if $n$ increases to a value that is no longer sufficiently small
relative to $d$. This is in line with
R.~A. Fisher's 1922 observation that
``more or less elaborate forms [of models]
will be suitable according to the volume of
data;'' cf.~\cite{Fis22a}.
And we stress that our results cannot guarantee that a given
simple model, like that discussed in the preceding paragraphs, is correct.
But we can guarantee, under the assumptions of Theorem~\ref{t0},
that most simple models are approximately correct,
in the sense that $\upsilon(B_d)$ is large and in the sense that
the left-hand sides of~\eqref{c1000} and~\eqref{c2000} are small,
provided only that $d$ is sufficiently large.
This also underscores the need for critical examination of the data and
of the model fit, irrespective of whether or not $d$ is large.
To this end, a very useful diagnostic tool is introduced by
Li in~\cite{Li97a},
namely a method to estimate, for a given unit-vector~$\beta$,
that unit-vector $\alpha$ for which the conditional mean of $\alpha'Z$
given $\beta'Z$ is most nonlinear in $\beta'Z$; see also
Section 6.1 in that reference.

The results obtained in this paper do not suggest that one should
abandon the search for complex and potentially nonlinear relations
in the data. But after such complex and/or nonlinear relations have
been accounted for, or in the case where none such can be found,
our results show how the use of simple linear models can be justified
without imposing strong regularity conditions on the true
data-generating process.

The discussion so far prompts for two extensions of our results that
are beyond the scope of this paper.
The first one is to extend our findings
to the case of more than one explanatory
variable, that is, the case where the conditioning is not on $\beta'
Z$ but
on $(\beta_1'Z, \beta_2'Z, \ldots, \beta_p' Z)$ for a collection of
$p$ mutually orthogonal unit-vectors $\beta_1,\ldots, \beta_p$. In fact,
Hall and Li~\cite{Hal93a} sketch an extension of their results to that
situation, so that an appropriate generalization of \eqref{c100} holds.
We will consider a corresponding generalization of
\eqref{c200} and also of the main results in \eqref{c1000} and~\eqref{c2000}
elsewhere.
The second extension is to provide explicit upper bounds for
the expressions on the left-hand sides of \eqref{c1000} and \eqref{c2000}
that converge to zero as $d\to\infty$ at a fast rate;
and to provide an explicit lower bound for $\upsilon(B_d)$
that converges to one as $d\to\infty$ also at a fast rate.

\subsection{Some consequences}
\label{S3}

\subsubsection{SIR, SAVE and related methods}
\label{ex2}

Many modern dimension reduction methods, like those based on inverse
conditional moments, rely on conditions like~(i) and~(ii) in
Section~\ref{S2}.
In particular, first-moment-based methods like
Sliced Inverse Regression~\cite{Li91a}
are based on a linear conditional mean requirement as in~(i).
(Besides, this requirement is also used in
several important results on generalized least squares under possible
link misspecification; see, for example,~\cite{Li89b} and the
references cited
therein.)
And second-moment-based techniques like
the Sliced Average Variance Estimator~\cite{Coo91a},
Principal Hessian Directions~\cite{Li92a}, or
Directional Regression~\cite{Li07b},
are based on both a linear conditional mean requirement as in (i),
and on a constant conditional variance requirement as in (ii).
Both conditions (i) and (ii) are also used in recent works such
as~\cite{Che10a,Chi02a,Don10a,Li12a}.

Given observations from a potentially rather complex data-generating process,
the dimension-reduction methods mentioned in the preceding paragraph
aim at finding a simpler
model that also describes the data. To justify the dimension reduction,
these methods
make assumptions to the effect that requirements like (i) or~(ii)
are satisfied, for \emph{one particular projection}, namely for the projection
on the so-called central subspace.
Under such assumptions, the central subspace or, equivalently,
the projection onto it, can be recovered from the data with good accuracy.
But as outlined in the \hyperref[sec1]{Introduction}, verifying such assumptions in practice
can be hard, particularly in situations where the sample size is
comparatively small.

Our results provide an alternative justification for requirements like
(i) and (ii). In particular,
in the setting of Theorem~\ref{t0}, we see that both (i) and (ii)
are approximately satisfied for \emph{most projections} $\beta' Z$
in the sense of \eqref{c1000} and \eqref{c2000},
provided only that the underlying dimension is large.
For the linear conditional mean condition, we stress that
the relation \eqref{c100} has been derived much earlier in~\cite{Hal93a}.

\subsubsection{Sparse linear modeling}
\label{ex1}

Consider the linear model with univariate response $y$ and
a $d$-vector of explanatory variables $w$, that is,
%
\begin{equation}
\label{m1} y = \theta' w + \varepsilon,
\end{equation}
where $\theta\in\R^d$ is unknown, and where
the error $\varepsilon$ has zero mean and constant variance conditional
on $w$.
We also assume that $y$ and $w$ are square integrable and centered so that
$\E w = 0$.
The leading case we have in mind is a situation where $d$, that is,
the number of available regressors, is as large as, or even much larger
than, the sample size.
To deal with such situations, it is common to assume that
the true model \eqref{m1} is equal to, or can be closely approximated by,
a ``sparse'' submodel that uses only a few explanatory
variables, and to use the available data to select and fit a sparse
submodel.
Such sparsity assumptions are clearly restrictive.
In the following, we argue that the results in this paper provide
weaker, and hence less restrictive, assumptions that also justify the fitting
of sparse submodels.

For illustration, consider now an extremely sparse model where $y$
is regressed on just one explanatory variable, say, $w_1$, that is,
%
\begin{equation}
\label{m2} y = c w_1 + e,
\end{equation}
where $c \in\R$ is unknown, and where $e$ has zero mean and constant
variance given~$w_1$.
To consider various possible justifications of the sparse
submodel \eqref{m2}, we first
rewrite the overall model \eqref{m1} as
\[
y = \bigl(\theta_1 w_1 + \E\bigl[
\theta_{\neg1}' w_{\neg1} \| w_1\bigr]
\bigr) + \bigl( \theta_{\neg1}' w_{\neg1} - \E\bigl[
\theta_{\neg1}' w_{\neg1} \| w_1\bigr] +
\varepsilon \bigr),
\]
where $\theta_{\neg_1}$ and $w_{\neg1}$ are obtained from
$\theta$ and $w$, respectively, by deleting the first component.

One possibility to justify the sparse model \eqref{m2} is to impose
the extreme
sparsity assumption
that all coefficients of $\theta_{\neg1}$ are zero, so that
$\theta_{\neg1}' w_{\neg1} = 0$. Then the relation
in the preceding display obviously reduces to $y = \theta_1 w_1 +
\varepsilon$,
and \eqref{m2} applies with $c = \theta_1$ and $e = \varepsilon$.
Under this extreme sparsity assumption we obtain, in particular, that the
sparse model \eqref{m2} is equivalent to the overall model \eqref{m1}
in terms of prediction, because $\E[y\|w_1] = \E[y \| w]$.
A slightly relaxed sparsity condition is to assume
that the coefficients of $\theta_{\neg1}$
are possibly nonzero but otherwise negligible,
that is, $\theta_{\neg1}' w_{\neg1} \approx0$, so that
the sparse model \eqref{m2} is approximately valid with
$c \approx\theta_1$ and $e \approx\varepsilon$.

An alternative justification of \eqref{m2} is to impose the assumption
that, given~$w_1$, the
conditional mean of $\theta_{\neg1}' w_{\neg1}$ is linear in $w_1$
and the conditional variance of
$\theta_{\neg1}' w_{\neg1}$ is constant in $w_1$.
In that case, the relation in the preceding display also reduces to
\eqref{m2}, but now with
$c= \Cov[ \theta' w , w_1] /\Var[w_1] =
\theta_1 + \sum_{i=2}^d \theta_i \Cov[w_1,w_i]/\Var[w_1]$,
and $e = \theta_{\neg1}' w_{\neg1} - \E[ \theta_{\neg1}' w_{\neg1}
\|w_1]
+ \varepsilon$.
Under this alternative assumption,
the model \eqref{m2} is valid but typically less accurate in terms of
prediction than the overall model \eqref{m1},
because, typically, $\E[y \| w_1] \neq\E[y\|w]$ and hence
$\Var[y\| w_1] > \Var[y\|w]$.
As before, these assumptions can
be relaxed by requiring that the conditional mean is approximately linear
and the conditional variance is approximately constant.

In the preceding two paragraphs,
we have considered two types of justifications for fitting the
submodel \eqref{m2}. Type (a): Exact or approximate sparsity assumptions.
And type (b): Exact or approximate linear conditional
mean and constant conditional variance assumptions. In practice, verifying
either of these assumptions for a given submodel can be
difficult. This raises the question as to which set of assumptions, that is,
(a) or (b), is more restrictive.
To this end, we first note that~(a) obviously implies (b).
For the more detailed comparison that we give in the following, we
assume that the law of $w$ is nondegenerate so that
$w$ can be written as $w = M Z$ for a $d$-vector $Z$ satisfying $\E Z =
0$ and
$\E Z Z' = I_d$. The $d\times d$ matrix~$M$ is a square root of the
variance/covariance matrix of $w$, which is nondegenerate by assumption;
which can be assumed to be symmetric;
and which need not be known in practice.
Then $\theta_{\neg1}' w_{\neg1}$ and $w_1$ can be written
as $\theta_{\neg1}' w_{\neg1} = a_\circ' Z$ and
$w_1 = b_\circ' Z$ with $a_\circ' = (0,\theta_2,\ldots, \theta_d) M$
and $b_\circ' = (1,0,\ldots, 0) M$.

The type (a) condition that $\theta_{\neg1} = 0$
entails that $a_\circ$
is equal to zero; the collection of $\theta$'s that satisfy this
condition
is the $1$-dimensional subset of the parameter space $\R^d$
that is spanned by $M^{-1} b_\circ= (1,0,\ldots,0)'$ (this collection depends
on~$b_\circ$).
More generally, for $y$ as in \eqref{m1} and for each vector $b$, the
simple model with response $y$ and with explanatory variable $b' Z$,
that is,
$y = c (b' Z) + e$, can be justified by the type (a) condition that
$\theta$ is parallel to $M^{-1}b$.
And, for each vector~$b$,
the approximate type (a) condition, that $\theta$ is approximately
parallel to $M^{-1} b$, is satisfied if $\theta$ belongs to
an appropriately small neighborhood of the span of $M^{-1} b$.

To study type (b) conditions on
the conditional moments of $\theta_{\neg1}' w_{\neg1}$
given $w_1$ or, equivalently,
on the conditional moments of $a_\circ'Z$ given $b_\circ'Z$,
we may replace the vectors $a_\circ$ and $b_\circ$ by standardized
versions $\alpha_\circ$ and $\beta_\circ$ that have length one,
for example, $\alpha_\circ= a_\circ/\|a_\circ\|$.
(Indeed, the conditional mean of $a_\circ'Z$
given $b_\circ'Z$ is linear, or approximately linear,
in $b_\circ'Z$, if and only if the same is true for the
conditional mean of $\alpha_\circ'Z$ given $\beta_\circ'Z$;
and a similar statement applies for the conditional variances, mutatis
mutandis.)
If Theorem~\ref{t0} applies and if $d$ is large,
then \emph{for most $\beta$'s and uniformly in $\alpha$},
the conditional mean of $\alpha'Z$ given $\beta' Z$
is approximately linear and the
conditional variance of $\alpha' Z$ given $\beta' Z$ is approximately
constant, in the sense of \eqref{c1000} and \eqref{c2000}.
In terms of the original parameter $\theta$, we note that uniformity in
$\alpha$ corresponds to uniformity in $\theta\in\R^d\setminus\{0\}$.

\section{Main result and outline of proof}
\label{S4}

Our main result is that \eqref{c100}--\eqref{c200} and also
\eqref{c1000}--\eqref{c2000} hold, for sets $B_d$ with
$\lim_{d\to\infty}\upsilon(B_d) = 1$.
We will establish this under the basic condition that
$Z$ has a Lebesgue density and is standardized so that
$\E Z = 0$ and $\E Z Z' = I_d$
for each $d$.
For the method of proof that we employ,
we also rely on two technical conditions,
which follow.

%
\begin{condition*}[(t1)]\hypertarget{co1}
For fixed $k\in\N$ and for each $d$, set $S_k = (Z_i' Z_j/d)_{i,j=1}^k$,
where the $Z_i$'s are i.i.d. copies of $Z$.
%
\begin{longlist}[(a)]
\item[(a)]
We have $\E[ (\sqrt{d}\Vert S_k-I_k\Vert )^{2 k + 1}] = O(1)$ as $d\to\infty$.
Moreover, let $H$ be a monomial in the elements of $S_k-I_k$ of degree
$h \leq2 k$.
If $H$ has a linear factor, then
$ d^{h/2} \E H = o(1)$.
And if $H$ consists only of quadratic factors in the
elements of $S_k-I_k$ above the diagonal, then
$ d^{h/2} \E H = 1+o(1)$.

\item[(b)]
Consider two monomials $G$ and $H$ of degree $g$ and $h$, respectively,
in the elements of $S_k -I_k$. If
$G$ is given by $Z_1'Z_2   Z_2'Z_3 \cdots Z_{g-1}'Z_g  Z_g'Z_1/d^g$,
if $H$ depends at least on those $Z_i$'s with $i\leq g$,
and if $2 \leq h < g \leq k$, then
$d^g \E G H = o(1)$.
\end{longlist}
\end{condition*}

\begin{condition*}[(t2)]\hypertarget{co2}
For fixed $k\in\N$, for each $d$,
and for any orthogonal $d\times d$
matrix $R$, the marginal density of the last $d-k+1$ components
of $R Z$ is bounded by ${d \choose k-1}^{1/2} B^{d-k+1}$ for
some constant $B$ that does not depend on $d$ or $R$.
\end{condition*}

Conditions~(\hyperlink{co1}{t1}) and~(\hyperlink{co2}{t2}) are always satisfied,
for any fixed $k$,
if the components of $Z$ are independent,
with bounded marginal densities and bounded marginal moments
of sufficiently high order; cf. Example~A.1 in the~supplementary material~\cite{Lee12b}.
Also, if $\E[(\sqrt{d} \Vert S_k - I_k\Vert )^{2 k+1}] = O(1)$, then
condition (\hyperlink{co1}{t1})(a) is satisfied if the elements of
$\sqrt{d}(S_k-I_k)$ jointly converge to a
Gaussian; cf. Example~A.2 in the~supplementary material~\cite{Lee12b}.
However, these conditions are more general than that and allow, in particular,
for situations where the components of $Z$ are dependent and/or
where the elements of $\sqrt{d}(S_k-I_k)$ do not converge in distribution.
Also note that both conditions are orthogonally invariant: If $Z$ satisfies
any one of them, then the same is true for any orthogonal
transformation of
$Z$.

The first requirement of condition~(\hyperlink{co1}{t1})(a) entails that
$d^{h/2} \E H = O(1)$
for any monomial $H$ in the elements of $S_k-I_k$ of degree up to $2 k
+ 1$.
Condition~(\hyperlink{co1}{t1})(b) strengthens parts of
condition~(\hyperlink{co1}{t1})(a) in the following sense:
Consider monomials $G$ and $H$
as in condition~(\hyperlink{co1}{t1})(b).
If condition~(\hyperlink{co1}{t1})(a) is satisfied,
then $\E G H = o(d^{-(g+h)/2})$ (because $G H$ is a monomial
of degree $g+h$ that has a linear factor).
Condition~(\hyperlink{co1}{t1})(b) requires that
$\E G H$ converges to zero at the faster rate $o(d^{-g})$.
Condition~(\hyperlink{co2}{t2}) ensures that the distribution of $Z$
is not ``too concentrated'' in certain directions, and is used together
with (\hyperlink{co1}{t1})(a) to guarantee uniform integrability of
$d/Z'Z$ and related quantities (see Proposition~E.1 in the~supplementary material~\cite{Lee12b}).
Also, our conditions should be compared to those used
in~\cite{Hal93a}.\footnote{\label{foot}
The results in~\cite{Hal93a} are stated under
high-level assumptions that are less specific but harder to verify;
see, for example,
conditions (3.21), (3.28) and (3.29) of Theorem 3.2 in that paper.
The only specific example that is actually shown in~\cite{Hal93a}
to satisfy
all three of these high-level conditions is the normal distribution;
cf. Example~4.2 and Remark~4.2 in that paper.}

We can now state the main result of this paper.

\begin{theorem}
\label{t0}
For each $d$, consider
a random $d$-vector $Z$ that has a Lebesgue density and that
is standardized such that $\E Z = 0$ and $\E Z Z' = I_d$.
If conditions~\textup{(\hyperlink{co1}{t1})(a)} and~\textup{(\hyperlink{co2}{t2})} are satisfied with $k=2$,
then there are Borel sets $B_d\subseteq\R^d$ satisfying
$\lim_{d\to\infty}\upsilon(B_d) = 1$, such that
\eqref{c1000} holds for each $\varepsilon>0$.
If condition~\textup{(\hyperlink{co1}{t1})} and condition~\textup{(\hyperlink{co2}{t2})} are
satisfied with $k=4$, then
the sets $B_d$ can be chosen so that also \eqref{c2000} holds for each
$\varepsilon>0$.
[Moreover, for each $x\in\R$ and each $\varepsilon>0$,
the relation \eqref{c100} obtains
under conditions~\textup{(\hyperlink{co1}{t1})(a)} and~\textup{(\hyperlink{co2}{t2})} with $k=2$,
and the relation \eqref{c200} holds
under conditions~\textup{(\hyperlink{co1}{t1})} and~\textup{(\hyperlink{co2}{t2})} with $k=4$.]
\end{theorem}

In the remainder of this section,
we give an outline of the proof of Theorem~\ref{t0}. The
proof is comprised of five main steps corresponding to the five propositions
that follow.

As the first step, it will be convenient to
replace the usual reference measure on~$\R^d$, that is, Lebesgue measure,
by the $d$-variate standard Gaussian measure, that is, $N(0,I_d)$.
The effect of this change of measure on conditional densities
and on conditional expectations involving $Z$ is described by the next result.

\begin{proposition}
\label{P1}
Fix $d\geq1$, and
consider a random $d$-vector $Z$ with Lebesgue density $f$.
Let $V \sim N(0, I_d)$, and write $\phi(\cdot)$ for the Lebesgue
density of
$V$. Moreover, for a fixed unit-vector $\beta\in\R^d$
and for each $x\in\R$, set
$W^{x|\beta} = x \beta+ (I_d-\beta\beta')V$.
Then the function $h(\cdot|\beta)$ defined by
\[
h(x|\beta) = \E \biggl[\frac{f( W^{x|\beta})}{\phi(
W^{x|\beta
})} \biggr]
\]
for $x\in\R$
is a density of $\beta'Z$ with respect to the univariate standard Gaussian
measure [i.e., $h(x|\beta) \phi_1(x)$ is a Lebesgue density of $\beta'Z$
if $\phi_1$ denotes the $N(0,1)$-density].
Moreover, if $\Psi\dvtx \R^d \to\R$ is such that $\Psi(Z)$ is
integrable, then a conditional expectation $\E[\Psi(Z)\|\beta'Z=x]$
of $\Psi(Z)$ given $\beta'Z$ satisfies
\[
\E \bigl[ \Psi(Z) \|\beta'Z =x \bigr] h(x|\beta) = \E \biggl[ \Psi
\bigl(W^{x|\beta}\bigr) \frac{f(W^{x|\beta})}{\phi(W^{x|\beta})} \biggr]
\]
whenever $x\in\R$ is such that $h(x|\beta) < \infty$.
\end{proposition}

This result allows us, for fixed $x\in\R$,
to study the marginal density of $\beta'Z$ at~$x$
as well as conditional expectations involving $Z$ given $\beta'Z=x$,
by considering unconditional means involving the random variable
$W^{x|\beta}$,
which has a $N(x\beta, I_d-\beta\beta')$-distribution.

Now, in order to derive~\eqref{c100}, we follow~\cite{Hal93a}
and use the following argument
(which can be traced back to Hoeffding~\cite{Hoe52a};
see also~\cite{Due11a}):
Set $\mu_{x|\beta} = \E[Z\| \beta'Z = x]$, and let $b$ be a random vector
in $\R^d$ with distribution $\upsilon$, that is, such that $b$
is uniformly distributed on the unit-sphere, that is, on the
set of unit-vectors in $\R^d$.
Then \eqref{c100} is equivalent to the statement that
$\Vert \mu_{x|b}\Vert^2 - x^2 $ converges to zero in probability
as a function of $b$, and this will follow if
\[
\E \bigl[ \bigl( \Vert \mu_{x|b}\Vert^2 -
x^2 \bigr) h^2(x|b) \bigr] \quad\mbox{and}\quad \E \bigl[ \bigl(
h(x|b) - 1 \bigr)^2 \bigr]
\]
both converge to zero as $d\to\infty$, where the expectations
are taken with respect to~$b$.
[Note that $\mu_{x|b}$ and $h(x|b)$ are measurable in view of
Corollary~B.2.
Also note that both integrands in the
preceding display are nonnegative.]
We now compute $\E[(h(x|b)-1)^2]$ as
%
\begin{eqnarray}\label{e1}
&&\int_{\beta} \biggl(\E \biggl[ \frac{f(W^{x|\beta})}{\phi(W^{x|\beta})}
\biggr] \biggr)^2 - 2 \E \biggl[ \frac{f(W^{x|\beta})}{\phi(W^{x|\beta})} \biggr] + 1
\upsilon(d \beta)
\nonumber
\\[-8pt]
\\[-8pt]
\nonumber
&&\qquad=  \E \biggl[ \frac{f(W_1)}{\phi(W_1)} \frac{f(W_2)}{\phi(W_2)} \biggr] - 2 \E
\biggl[ \frac{f(W_1)}{\phi(W_1)} \biggr] + 1,
\end{eqnarray}
where the $W_i$'s are defined as $W_i = x b + (I_d - b b' )V_i$, $i=1,2$,
with $V_1$ and $V_2$ i.i.d. $N(0,I_d)$ and independent of~$b$.
Note that
$W_1$ and $W_2$ are dependent because both share
the same random unit-vector~$b$.
If \eqref{e1} converges to zero or, equivalently, if $h(x|b) \to1$
in $L^2(b)$ as $d\to\infty$, then $h(x|b) < \infty$ with probability one
for sufficiently large $d$. For such $d$,
we see that the second statement of Proposition~\ref{P1} applies
for
$\upsilon$-almost all $\beta$, such that
$\E[(\Vert \mu_{x|b}\Vert^2 - x^2) h^2(x|b)]$ can be written as
%
\begin{equation}
\label{e2} \E \biggl[ \bigl(W_1' W_2 -
x^2 \bigr) \frac{f(W_1)}{\phi(W_1)} \frac{f(W_2)}{\phi(W_2)} \biggr],
\end{equation}
by arguing as in the derivation of \eqref{e1},
as is easy to see.
With this, we see that \eqref{c100} holds if both
\eqref{e1} and \eqref{e2} go to zero as $d\to\infty$.
And if \eqref{e1} and \eqref{e2} converge to zero uniformly in $x$ over
compact subsets of $\R$, that is, if
the suprema of \eqref{e1} and \eqref{e2} over $x$ satisfying
$|x| \leq M$ converge to zero as $d\to\infty$ for each $M>0$,
then it is not difficult to also derive \eqref{c1000}
by employing standard arguments; for details,
see Lemma~B.4(i) in the~supplementary material~\cite{Lee12b}.

To establish \eqref{c200}, we employ a similar strategy and
write $\Delta_{x|\beta}$ as shorthand for
the $d\times d$ matrix
$\Delta_{x|\beta} = \E[Z Z'\| \beta'Z = x] - (I_d + (x^2-1)\beta
\beta')$.
With $b$ again uniform on the unit-sphere,
the goal is to show convergence of
the largest singular value $\Vert \Delta_{x|b}\Vert $
to zero in probability with respect to $b$.
But this follows if
\[
\trace\Delta^{k}_{x|b} \stackrel{p} {\longrightarrow} 0
\]
as $d\to\infty$ for some even integer $k$ (where, again,
Corollary~B.2 guarantees measurability).
This, and hence also \eqref{c200}, will follow if
\[
\E \bigl[ \trace\Delta^k_{x|b}
h^k(x|b) \bigr]
\]
converges to zero as $d\to\infty$ for some even integer $k$
and if, in addition, also \eqref{e1} converges to zero.
(We shall find, at the end of the section, that the expression
in the preceding display converges to zero for $k=4$ but typically
not for $k=2$.)
To analyze the expectation in the preceding display, define the function
$\Delta_{x|\beta}(z)$ as
$\Delta_{x|\beta}(z) = z z' - (I_d+(x^2-1)\beta\beta')$ for $z\in
\R^d$.
We now argue as in the preceding paragraph:
Assume that \eqref{e1} converges to zero, and
assume that $d$ is sufficiently large so that $h(x|\beta)<\infty$
for $\upsilon$-almost all $\beta$. For such $d$,
use Proposition~\ref{P1} to see that $\E[\trace\Delta^k_{x|b}
h^k(x|b)]$ equals
\[
\E \biggl[ \trace\Delta_{x|b}(W_1)\cdots
\Delta_{x|b}(W_k) \frac{f(W_1)}{\phi(W_1)} \cdots \frac{f(W_k)}{\phi(W_k)}
\biggr],
\]
where, similarly to before,
$W_i = x b + (I_d-b b') V_i$ with the $V_i$, $1\leq i \leq k$, i.i.d.
$N(0,I_d)$ independent of $b$.
Instead of computing the trace in the preceding display directly,
we find it convenient to break it into smaller, more manageable, pieces.
Indeed, we find that the expression in the preceding display can be
written as the
weighted sum of the terms
%
\begin{eqnarray}
\label{e3} \qquad &\displaystyle\sum_{j=1}^{k} \pmatrix{k \cr j}
(-1)^{j} \E \biggl[ \bigl(W_1'
W_2 \cdots W_{j}' W_1 -
d+1-x^{2 j} \bigr) \frac{f(W_1)}{\phi(W_1)} \cdots \frac{f(W_k)}{\phi(W_k)}
\biggr],&
\nonumber
\\[-8pt]
\\[-8pt]
\nonumber
&\displaystyle  \E \Biggl[ \Biggl( \prod_{i=1}^m
W_{j_{i-1}+1}' W_{j_{i-1}+2} \cdots W_{j_{i}-1}'
W_{j_{i}} - x^{2(j_m - m)} \Biggr) \frac{f(W_1)}{\phi(W_1)} \cdots
\frac{f(W_k)}{\phi(W_k)} \Biggr]&
\end{eqnarray}
for $m\geq1$ and indices $j_0,\ldots, j_m$ satisfying
$j_0=0$, $j_m < k$, and
$j_{i-1}+1 < j_{i}$ whenever $0 < i \leq m$.
In addition, we find that the weights in this weighted sum
depend only on $k$ and on $x$, and that the weights
are continuous in $x\in\R$; cf. Lemma~B.3 for details.
[Note that, in \eqref{e3},
we write $W_1' W_2 \cdots W_j'W_1$ as shorthand
for $\trace\prod_{i=1}^j W_i W_i'$, and we also use notation like
$W_1'W_2\cdots W_{j-1}' W_j$ as shorthand
for $\prod_{i=1}^{j-1} W_i' W_{i+1}$.]
Hence, \eqref{c200} holds if the expression in \eqref{e1} and
those in \eqref{e3}
go to zero as $d\to\infty$, the latter for some even integer $k$,
and for any $m$ and $j_0,\ldots, j_{m}$ as indicated.
Moreover, \eqref{c2000} holds provided that
the expressions in \eqref{e1} and \eqref{e3} all converge to zero uniformly
in $x$ over compact subsets of $\R$; details are given in
Lemma~B.4(ii) in the~supplementary material~\cite{Lee12b}.

To understand the large-$d$-behavior of the quantities in
\eqref{e1}, \eqref{e2}, and \eqref{e3},
we need to understand the joint distribution of the $W_j$'s,
which is described by the next result.

\begin{proposition}
\label{P2}
For $d$ and $k$ satisfying $1 \leq k < d$,
the joint distribution of $W_1,\ldots, W_k$
has a density with respect to Lebesgue measure that we denote
by $\varphi_x(w_1,\ldots, w_k)$, and this density satisfies
\[
\frac{
\varphi_x(w_1,\ldots,w_k)
}{
\phi(w_1) \cdots\phi(w_k)
}  = \frac{(d/2)^{-k/2}\Gamma(d/2)}{\Gamma((d-k)/2)} \det S_k^{-1/2}
\biggl( 1 - \frac{x^2}{d} \iota' S_k^{-1}
\iota \biggr)^{{(d-k-2)}/ {2}} e^{({k}/{2}) x^2}
\]
if $S_k$ is invertible with $x^2 \iota' S_k^{-1} \iota< d$,
and $\varphi_x(w_1,\ldots, w_k) = 0$ otherwise, where
$S_k = (w_i'w_j/d)_{i,j=1}^k$ denotes the $k\times k$ matrix
of scaled inner products of the $w_i$'s, and
$\iota=(1,\ldots,1)'$ denotes an appropriate vector of ones.
\end{proposition}

Using Proposition~\ref{P2}, we can rewrite the quantities of
interest in \eqref{e1}, \eqref{e2}, and \eqref{e3} as follows:
For example, \eqref{e2}
equals
%
\begin{eqnarray}\label{a}
&& \int\!\!\!\int\bigl(w_1' w_2 -
x^2\bigr) \frac{f(w_1)}{\phi(w_1)} \frac{f(w_2)}{\phi(w_2)} \varphi_x(w_1,w_2)
\,d w_1 \,d w_2
\nonumber\\
&&\qquad
= \int\!\!\!\int\bigl(w_1' w_2 -
x^2\bigr) \frac{\varphi_x(w_1,w_2)}{\phi(w_1) \phi(w_2)} f(w_1) f(w_2)
\,d w_1 \,d w_2
\\
&&\qquad = \E \biggl[ \bigl(Z_1' Z_2 -
x^2 \bigr) \frac{\varphi_x(Z_1,Z_2)}{\phi(Z_1) \phi(Z_2)}
\biggr],\nonumber
\end{eqnarray}
where $Z_1$ and $Z_2$ are i.i.d. copies of $Z$. In a similar fashion,
the quantities in \eqref{e3} can be written as
%
\begin{eqnarray}
\label{b}
&\displaystyle\sum_{j=1}^{k} {k \choose j}
(-1)^j \E \biggl[ \bigl(Z_1'
Z_2 \cdots Z_{j}' Z_1 -
d+1-x^{2 j} \bigr) \frac{\varphi_x(Z_1,\ldots, Z_k)}{ \phi(Z_1)\cdots\phi(Z_k)}
\biggr],&
\nonumber
\\[-8pt]
\\[-8pt]
\nonumber
&\displaystyle \E \Biggl[ \Biggl( \prod_{i=1}^m
Z_{j_{i-1}+1}' Z_{j_{i-1}+2} \cdots Z_{j_{i}-1}'
Z_{j_{i}} - x^{2(j_m - m)} \Biggr) \frac{\varphi_x(Z_1,\ldots, Z_k)}{ \phi(Z_1)\cdots\phi(Z_k)}
\Biggr]&
\end{eqnarray}
for $Z_i$, $i=1,\ldots, k$, i.i.d. as $Z$. And finally
\eqref{e1} reduces to
%
\begin{equation}
\label{c} \qquad \E \biggl[ \frac{\varphi_x(Z_1, Z_2)}{\phi(Z_1)\phi(Z_2)} \biggr] - 2 \E \biggl[ \frac{\varphi_x(Z_1)}{\phi(Z_1)}
\biggr] + 1.
\end{equation}

Recall that our goal is to show that the expressions in \eqref
{a}--\eqref{c}
converge to zero as $d\to\infty$, uniformly in $x$ over compact
subsets of
$\R$. We show, in fact, a slightly stronger statement, motivated by the
obvious
conjecture that the expected value of density ratios in \eqref
{a}--\eqref{c},
like $\varphi_x (Z_1,Z_2)/(\phi(Z_1)\phi(Z_2))$, converges to one,
and also by the observation that the expression in \eqref{a} is a
special case
of the second expression in \eqref{b} with $k$ replaced by $2$.
For an even integer $k$ and
for each $l=1,\ldots, k$, for each $m\geq0$ and for any
indices $j_0,\ldots, j_m$
that satisfy\vadjust{\goodbreak} $j_0=0$, $j_m\leq l$, and $j_{i-1}+1 < j_i$ whenever
$0<i\leq m$,
consider the expressions
%
\begin{eqnarray}\qquad
\label{B} \E \Biggl[ \Biggl( \prod_{i=1}^m
Z_{j_{i-1}+1}' Z_{j_{i-1}+2} \cdots Z_{j_{i}-1}'
Z_{j_{i}} \Biggr) \frac{\varphi_x(Z_1,\ldots, Z_l)}{ \phi(Z_1)\cdots\phi(Z_l)} \Biggr] - x^{2(j_m-m)}
\end{eqnarray}
and also the expression
%
\begin{eqnarray}\qquad
\label{C} &\sum_{j=1}^{k} \pmatrix{k \cr j}
(-1)^j \E \biggl[ \bigl(Z_1'
Z_2 \cdots Z_{j}' Z_1 -d \bigr)
\frac{\varphi_x(Z_1,\ldots, Z_k)}{ \phi(Z_1)\cdots\phi(Z_k)} \biggr] - \bigl(1-x^2\bigr)^k .
\end{eqnarray}
Convergence to zero of the expressions in \eqref{B} corresponding to
$k = 2$ (uniformly over
compacts in $x$) entails convergence to zero of \eqref{a} and \eqref{c}
(uniformly over compacts in $x$).
And if both \eqref{B} and \eqref{C} converge to zero for some even
integer~$k$
(uniformly over
compacts), then the expressions in \eqref{b} corresponding to that
same $k$
also converge to zero (uniformly over compacts).
[Convergence to zero of \eqref{a} follows
from convergence to zero of \eqref{B} with
$m=1$, $j_m=2$ and $l=2$ together with convergence to zero of \eqref{B}
with $m=0$ and $l=2$.
Convergence to zero of~\eqref{c} follows from
convergence to zero of~\eqref{B} with $m=0$ and $l=1$ together
with convergence to zero of~\eqref{B} with $m=0$ and $l=2$.
For the first expression in \eqref{b}, convergence to zero
follows from convergence to zero of~\eqref{B} with $m=0$ and $l=k$
and from convergence to zero of \eqref{C}, in view of the binomial
theorem. Similarly, for the second expression in \eqref{b},
convergence to zero follows from convergence to zero of \eqref{B}
and of convergence to zero of the special case
of \eqref{B} where $m=0$ and $l=k$.]

The expressions in \eqref{B}--\eqref{C} both involve expected values of
a polynomial in $Z_i' Z_j$ for some pairs $(i,j)$,
multiplied by $\varphi_x(Z_1,\ldots, Z_l) / (\phi(Z_1)\cdots\phi(Z_l))$,
$l=1,\ldots, k$. To proceed, we need the polynomial approximation to
$\varphi_x(Z_1,\ldots, \break Z_l) /  (\phi(Z_1)\cdots\phi(Z_l))$
that is provided by the next result.

\begin{proposition}
\label{P3}
Fix $M>0$ and $x$ satisfying $|x|\leq M$. Moreover,
consider integers $k$ and $d$ that satisfy
$k\geq1$ and
$d > \max\{ 3 k, 2 (k+1) M^2\}$,
and $d$-vectors $w_1,\ldots, w_k$ that are such that the $k\times k$ matrix
$S_k = (w_i'w_j/d)_{i,j=1}^k$ satisfies $\Vert S_k - I_k\Vert  < 1/(2 k)$. Then
$\varphi_x(w_1,\ldots, w_k)$ is such that
\[
\frac{\varphi_x(w_1,\ldots, w_k)}{ \phi(w_1) \cdots\phi(w_k)} = \psi_x(S_k - I_k) +
\Delta,
\]
where $\psi_x(S_k-I_k)$ is a polynomial of degree up to $k$ in the
elements of $S_k-I_k$.
The coefficients of the polynomial $\psi_x(\cdot)$ depend on
$k$, $x$ and $d$, and are bounded,
in absolute value and uniformly in $x \in[-M,M]$,
by a constant $C$, and
$\Delta$ satisfies $|\Delta| < D \Vert S_k - I_k\Vert^{k+1}$, for some
constants $C = C(k,M)$ and $D=D(k,M)$ that depend
only on $k$ and on $M$. Moreover, both $\psi_x(S_k-I_k)$ and $\Delta$
are invariant under permutations of the $w_i$'s so that, in particular,
$\psi_x(S_k-I_k)$ is unchanged when $S_k$ is replaced by
the matrix $(w_{\pi(i)}' w_{\pi(j)}/d)_{i,j=1}^k$\vadjust{\goodbreak} for any permutation
$\pi$ of $k$ elements.
[The coefficients of $\psi_x(\cdot)$ and the bounds $C$ and $D$
can be obtained explicitly upon inspection of the proof.]
\end{proposition}

When studying the expected values in \eqref{B}--\eqref{C},
Proposition~\ref{P3} suggests that the density ratio
$\varphi_x(Z_1,\ldots, Z_k)/(\phi(Z_1)\cdots\phi(Z_k))$
can be approximated by the polynomial
$\psi_x(S_k-I_k)$.
The resulting approximations to \eqref{B} and \eqref{C} are
%
\begin{equation}
\label{B1} \E \Biggl[ \Biggl( \prod_{i=1}^m
Z_{j_{i-1}+1}' Z_{j_{i-1}+2} \cdots Z_{j_{i}-1}'
Z_{j_{i}} \Biggr) \psi_x(S_l-I_l)
\Biggr] - x^{2(j_m-m)}
\end{equation}
and
\begin{equation}\quad
\label{C1} \sum_{j=1}^{k} \pmatrix{k \cr j}
(-1)^j \E \bigl[ \bigl( Z_1'
Z_2 \cdots Z_{j}' Z_1 - d \bigr)
\psi_x(S_k-I_k) \bigr] -
\bigl(1-x^2\bigr)^k,
\end{equation}
respectively.
For these approximations to be useful, we need to show that the difference
of \eqref{B} and \eqref{B1},
and also the difference of \eqref{C} and \eqref{C1},
converges to zero as $d\to\infty$, uniformly in $x$ over compact subsets
of $\R$. The technical difficulty here is that
expressions like, for example, $Z_1' Z_2 \cdots Z_j' Z_1 - d$ in \eqref{C1}
have zero mean but do not converge to zero in probability.
To deal with this, we
rely on conditions~(\hyperlink{co1}{t1})(a) and~(\hyperlink{co2}{t2}).

\begin{proposition}\label{P4}
For each $d$,
consider a random $d$-vector $Z$ that has a Lebesgue density, that
is standardized such that $\E Z = 0$ and $\E Z Z' = I_d$, and
that satisfies conditions~\textup{(\hyperlink{co1}{t1})(a)} and~\textup{(\hyperlink{co2}{t2})} for
some fixed integer $k$.
Let $H(S_k-I_k)$ be a (fixed) monomial in the elements of
$S_k-I_k$ whose degree, denoted by $\deg(H)$,
satisfies $0\leq\deg(H) \leq k$. Then
\[
\E \biggl[ d^{{(k+\deg(H))}/ {2}}\bigl | H(S_k-I_k)\bigr|
\biggl\llvert \frac{\varphi_x(Z_1,\ldots, Z_k)}{\phi(Z_1)\cdots\phi(Z_k)} - \psi_x(S_k -
I_k)\biggr\rrvert \biggr]
\]
converges to zero as $d\to\infty$, uniformly in $x$ over compact subsets
of $\R$.
\end{proposition}

If Proposition~\ref{P4} applies, then it is not difficult to see
that the difference between \eqref{B} and \eqref{B1},
and also the difference between \eqref{C} and \eqref{C1},
converges to zero, uniformly in $x$ over compact subsets of $\R$,
as required.
For example, consider the difference of \eqref{C} and \eqref{C1},
which both involve $k$ expected values indexed by $j=1,\ldots, k$,
and focus on the
difference of those expected values corresponding to the index $j$.
Also, recall that $k$ is an even integer, so that $k>1$.
For $j=1$, we simply use Proposition~\ref{P4} with the
monomial $(S_k-I_k)_{1,1}$, and note that
$|Z_1'Z_1-d| = d| (S_k-I_k)_{1,1}| < d^{(k+1)/2}| (S_k-I_k)_{1,1}|$.
For $j>1$, we first write
$|Z_1'Z_2 \cdots Z_j' Z_1|$ as
\[
d^j \bigl\llvert (S_k-I_k)_{1,2}
\cdots(S_k-I_k)_{j,1}\bigr\rrvert \leq
d^{(j+k)/2} \bigl\llvert (S_k-I_k)_{1,2}
\cdots(S_k-I_k)_{j,1} \bigr\rrvert .
\]
Now use Proposition~\ref{P4} twice, first with the
monomial $(S_k-I_k)_{1,2} \cdots(S_k-I_k)_{j,1}$ of degree $j\leq k$,
and then\vadjust{\goodbreak} with the degree-zero monomial, and note that
$d \leq d^{(k+0)/2}$ here,
to see that the difference of expected values corresponding to the
index $j>1$
also converges to zero, uniformly in $x$ over compact subsets of $\R$.
The difference of \eqref{B} and \eqref{B1} is treated similarly.

To show that the expressions in \eqref{B1} and \eqref{C1} converge to zero,
the following observation will be useful: If the $Z_i$'s in \eqref{B1} and
\eqref{C1} are replaced by independent standard normal vectors $V_i$
(and if
$S_l$ and $S_k$ are replaced by the corresponding Gram matrices of the $V_i$'s),
then the resulting expressions both converge to zero as $d\to\infty$,
uniformly in $x$ over compact subsets of $\R$.
To establish convergence to zero of \eqref{B1}--\eqref{C1}, uniformly on
compacts in $x$, it
therefore is sufficient to study the differences
between the expressions in \eqref{B1}--\eqref{C1}
and the same expressions with the $Z_i$'s
replaced by $V_i$'s that are i.i.d. standard normal,
and to show that these differences converge to zero as $d\to\infty$,
uniformly in $x$ over compacts subsets of $\R$.
(To derive the last observation, we first note that \eqref{B} and
\eqref{C},
with the $Z_i$'s replaced by the $V_i$'s, are both equal to zero.
Indeed, with this replacement, the expectation in \eqref{B} is equal to
$\E[ \prod_{i=1}^m
W_{j_{i-1}+1}' W_{j_{i-1}+2} \cdots W_{j_{i}-1}' W_{j_{i}}
]$,
because $\phi(v_1)\cdots\phi(v_l)$ is the joint density of
$V_1,\ldots,
V_l$, and
because $\varphi_x(w_1,\ldots, w_l)$ is the joint density of
$W_1,\ldots, W_l$. Conditional on $b$,
the $W_i$'s are conditionally i.i.d.,
with $\E[W_i \| b] = x b$ and
$\E[W_i W_i'\|b] = I_d + (x^2 -1) b b'$. In view of this, it is elementary
to verify that \eqref{B} with the $Z_i$'s replaced by the $V_i$'s
is equal to zero. A similar argument applies, mutatis mutandis,
to \eqref{C}.
Next, we note that Proposition~\ref{P4} applies if $Z$ is replaced by
a standard normal vector $V$ [that conditions~(\hyperlink{co1}{t1})(a)
and~(\hyperlink{co2}{t2})
hold when $Z$ is replaced by $V$ follows either from Example~A.1
or upon a simple direct computation]. When replacing $Z$ by $V$ throughout,
this entails that \eqref{B1} and~\eqref{C1} converge to the same limit
as \eqref{B} and \eqref{C}, that is, to zero, uniformly over compacts in
$x$.)

To put this idea to work, expand $\psi_x(S_k-I_k)$ into a weighted sum of
monomials (in the elements of $S_k-I_k$), where the weight of each monomial
in the sum is given by the coefficient of that monomial in
$\psi_x(S_k-I_k)$; similarly, $\psi_x(S_l-I_l)$ can also be written
as a weighted sum of such monomials for each $l\leq k$.
We see that the integrand in \eqref{B1} for $m=0$, that is,
$\psi_x(S_l-I_l)$, can be written as the weighted sum of
monomials in the elements of $S_k - I_k$. Similarly,
for $Z$ as in Theorem~\ref{t0},
the integrand in \eqref{B1} for $m>0$
can be written as the weighted sum of expressions of the form
%
\begin{equation}
\label{G0} d^{\deg(G)} \bigl( G - \E[G] \bigr) H
\end{equation}
for two monomials $G$ and $H$ in the elements of $S_k-I_k$ of degree $k$
or less,
where~$G$ is given by the monomial
%
\begin{equation}
\label{openc} \prod_{i=1}^m
(S_k-I_k)_{j_{i-1}+1,j_{i-1}+2} \cdots(S_k-I_k)_{j_i-1,j_i}
\end{equation}
of degree $j_m - m$,
for some $m > 0$ and indices $j_0,\ldots, j_m$ that satisfy
$j_0=0$, $j_m \leq k$, and $j_{i-1} +1<j_i$ whenever\vadjust{\goodbreak} $0<i\leq m$.
In this weighted expansion of~\eqref{B1}, note that the weight of
\eqref{G0} depends on
$x$, on $H$ (through its degrees) and also on~$d$, in such a way that
the weight is bounded in $x$ and $d$ as long as $x$ is restricted to a
compact set (cf. Proposition~\ref{P3}).
Lastly, consider the integrand in \eqref{C1}, for $Z$ as in
Theorem~\ref{t0}.
Arguing as before, we can
write that integrand as the weighted sum of expression of the form
\eqref{G0}, where here $G$ is given by the monomial
%
\begin{equation}
\label{closedc} (S_k-I_k)_{1,2}
\cdots(S_k-I_k)_{j-1,j} (S_k-I_k)_{j,1}
\end{equation}
of degree $j$
for some $j$ satisfying $1\leq j \leq k$.
And, again, in this weighted sum, the weight of each term depends
on $x$ and on $H$ (through its degrees), and that weight is
bounded in $x$ and $d$ over compacts in $x$.
(Note that, under the assumptions of Theorem~\ref{t0}, it is elementary
to verify that $\E[G] = 0$ whenever $G$ is given by~\eqref{openc} and also
whenever $G$ is given by~\eqref{closedc} with
$j=1$, and that $d^{\deg(G)}\E[G] = d$ if $G$ is
given by~\eqref{closedc} with $j>1$.)

\begin{proposition}\label{P5}
For each $d\geq1$, assume that $Z$ is as in Theorem~\ref{t0}.
Fix an integer $k\geq1$, and
let $G$ and $H$ be two (fixed)
monomials in the elements of $S_k - I_k$ of degree $k$ or less,
define $G^\ast$ and $H^\ast$ as $G$ and $H$,
respectively, but with the $Z_1,\ldots, Z_k$ replaced by i.i.d. standard
Gaussian $d$-vectors, and consider
%
\begin{equation}
\label{p51} \E \bigl[ d^{\deg(G)} \bigl( G - \E[G] \bigr) H \bigr] - \E
\bigl[ d^{\deg(G^\ast)} \bigl( G^\ast- \E\bigl[G^\ast\bigr]
\bigr) H^\ast \bigr].\vspace*{-9pt}
\end{equation}
%
%
\begin{longlist}[(ii)]
\item[(i)]
Assume that condition~\textup{(\hyperlink{co1}{t1})(a)} applies with
the integer $k$ as chosen here, and that $k\leq4$.
Then $\E[H] - \E[H^\ast]$ and also the expression in \eqref{p51}
converge to zero as $d\to\infty$ for each monomial $G$ as in
\eqref{openc}.
\item[(ii)]\label{p5ii}
Assume that condition~\textup{(\hyperlink{co1}{t1})} is
satisfied with the integer $k$ as chosen here.
Let~$G$ be given by the monomial in \eqref{closedc}
for some $j$, $1\leq j \leq k$.
Then the expression in
\eqref{p51} converges to zero as $d\to\infty$, unless
either \textup{(a)} $H=(S_k-I_k)_{a,a}$ for some $a$ satisfying $1\leq a \leq j$,
\textup{(b)} $H=(S_k - I_k)_{a,b}$ with $1\leq a < b \leq j$, or~\textup{(c)}
$H=( (S_k - I_k)_{a,b})^2$ with $1 \leq a < b \leq j$.
In case \textup{(a)}, the expression in \eqref{p51} is
equal to $\Var[Z_1'Z_1]/d - 2$; in case \textup{(b)}, it is equal
to $\E[ (Z_1'Z_2)^3]/d$; and in case~\textup{(c)}, it equals
$\Var[(Z_1'Z_2)^2]/d^2 - 2(1+3/d)$.
\end{longlist}
\end{proposition}

To complete the proof of Theorem~\ref{t0}, let $Z$ be as in the theorem.
We first assume that conditions~(\hyperlink{co1}{t1})(a) and~(\hyperlink{co2}{t2})
are satisfied with $k=2$.
The relation \eqref{c100} holds for each $x$ and $\varepsilon$,
if we can show that the expression in \eqref{B1}
converges to zero for each collection of
indices $l$, $m$, $j_0,\ldots, j_m$
so that $1\leq l\leq2$, $m\geq0$, $j_0=0$, $j_m \leq l$, and
$j_{i-1}+1 < j_i$ for each $i=1,\ldots, m$.
Moreover, \eqref{c1000} holds for each $\varepsilon>0$, if
convergence zero of \eqref{B1} is uniform in $x$ over compacts.
To this end,
consider the difference of the expression in \eqref{B1} and
of the same expression with the $Z_i$'s replaced by $V_i$'s that are i.i.d.
standard normal $d$-vectors. Expanding the polynomial $\psi_x$ into a
weighted sum of monomials, the difference in question can be written
as a weighted sum of expressions of the form
$\E[H] - \E[H^\ast]$ in case \mbox{$m=0$}, and
as a weighted sum of expressions of the form
\eqref{p51} in case \mbox{$m>0$}, where the
weight is given by the coefficient of the monomial $H$ in $\psi_x$,
and where $G$ is of the form \eqref{openc}
[the monomials $H$, $G^\ast$ and $H^\ast$ are as in
Proposition~\ref{P5}(i)].
By Proposition~\ref{P3}, we see that the coefficients of $\psi_x$
are bounded uniformly in $d$ and uniformly in $x$ over compacts. And by
Proposition~\ref{P5}(i),
we see that expressions of the form $\E[H]-\E[H^\ast]$
or of the form \eqref{p51} with $G$ as in \eqref{openc} all converge to
zero. Therefore, \eqref{B1} converges to zero, uniformly in $x$ over
compacts subsets of $\R$.

Finally, we assume that conditions~(\hyperlink{co1}{t1}) and~(\hyperlink{co2}{t2})
hold with $k=4$. To derive \eqref{c200}
for fixed $x$ and $\varepsilon$, we show that the expressions
in \eqref{B1} and \eqref{C1} converge to zero [with the
indices $l$, $m$, $j_0,\ldots, j_m$ in \eqref{B1} now
so that $1\leq l\leq4$, $m\geq0$, $j_0=0$, $j_m \leq l$, and
$j_{i-1}+1 < j_i$ for each $i=1,\ldots, m$].
And \eqref{c2000} holds for each $\varepsilon>0$ if convergence in
\eqref{B1} and \eqref{C1} is uniform in $x$ over compact sets.
Now convergence to zero of \eqref{B1} (with $k=4$ here), uniformly over
compacts, follows by arguing as in the preceding paragraph, mutatis mutandis.
To deal with \eqref{C1},
consider the difference of \eqref{C1} and of the same expression
with the $Z_i$'s replaced by i.i.d. standard Gaussian $V_i$'s.
Again, this can be written as a weighted sum of expressions of the
form \eqref{p51}, with $G$ now as in \eqref{closedc},
where the weights are bounded uniformly in $d$ and uniformly in $x$
over compacts in view of Proposition~\ref{P3}. And by
Proposition~\ref{P5}(ii),
we see that \eqref{p51} converges to zero except for those
$H$'s that correspond to the cases (a), (b) and (c) in
Proposition~\ref{P5}(ii).
Write ${\cal H}_a$, ${\cal H}_b$, and ${\cal H}_c$ for the
collection of all monomials $H$ where the case (a), (b), or (c)
of Proposition~\ref{P5}(ii) occurs, respectively.
For each $H \in{\cal H}_a$, the value of \eqref{p51} is
given $\Var[Z_1'Z_1]/d - 2$ and hence does not depend on $H$
in view of Proposition~\ref{P5}(ii).
And for each $H \in{\cal H}_a$, the coefficient of $H$ in the
polynomial $\psi_x(S_k - I_k)$ also does not depend on $H$
in view of Proposition~\ref{P3}, because the monomials in ${\cal H}_a$
can be obtained from each other by permutations (or re-labelings)
of $Z_i$'s.
Consider now the difference of \eqref{C1} and the same expression with
the $Z_i$'s replaced by i.i.d. standard Gaussians.
The combined contribution of the monomials in ${\cal H}_a$ to that
difference is given by
\[
\sum_{j=1}^k \pmatrix{k \cr j}
(-1)^j j = -k \sum_{j=0}^{k-1}
\pmatrix{k-1 \cr j} (-1)^j
\]
multiplied by a constant (namely by $\Var[Z_1'Z_1]/d - 2$ times
the common coefficient of the monomials from
${\cal H}_a$ in $\psi_x(S_k-I_k)$).
Similarly, the monomials in ${\cal H}_b$ contribute
\[
\sum_{j=1}^k \pmatrix{k \cr j}
(-1)^j \pmatrix{j \cr 2} = \sum_{j=0}^{k-2}
\pmatrix{k-2 \cr j} (-1)^j
\]
multiplied by a constant. And the combined contribution of the
monomials in ${\cal H}_c$ is also given by the expression in the
preceding display multiplied by another constant.
Because we have $k=4$ here, the expressions in the last two displays
are both equal to zero.
Except for the more technical arguments that we collect in the~supplementary material~\cite{Lee12b},
this concludes the proof of Theorem~\ref{t0}.

%
%


\begin{supplement}[id=suppA]
\stitle{Appendix: Proofs for Section~\ref{S4}}
\slink[doi]{10.1214/12-AOS1081SUPP} 
\sdatatype{.pdf}
\sfilename{aos1081\_supp.pdf}
\sdescription{The Appendix contains several more technical arguments
that are used in Section~\ref{S4} including, in particular,
Examples~A.1 and~A.2, as well as
the proofs of Propositions~\ref{P1} through~\ref{P5}.}
\end{supplement}

%

\printaddresses

\end{document}